\documentclass[12pt]{article}
\usepackage{epsfig}
\usepackage{graphicx}
\usepackage{amsmath}
\usepackage{amsfonts}
\usepackage{amssymb}
\def\theequation{\arabic{section}.\arabic{equation}}
\newtheorem{theorem}{Theorem}[section]

\newtheorem{lemma}[theorem]{Lemma}

\newtheorem{proposition}{Proposition}

\newcommand\hE{\hat{E}}

\newcommand\hJ{{\hat J}}

\newcommand\BbR{\mathbb{R}}

\newcommand{\be}{\begin{equation}}
\newcommand{\ee}{\end{equation}}
\newcommand{\bea}{\begin{eqnarray}}
\newcommand{\eea}{\end{eqnarray}}
\newcommand{\beaa}{\begin{eqnarray*}}
\newcommand{\eeaa}{\end{eqnarray*}}

\newcommand\ts{\tilde{s}}
\newcommand\veps{\varepsilon}

\newcommand\p{{\partial}}
\newcommand\di{{\rm d}}

\begin{document}

\author{
Wenjun Zhang\footnote{School of Computing and Mathematical Sciences,
Auckland University of Technology,
Private Bag 92006, Auckland 1142, New Zealand.
Phone: +64-9-9219999-x5094; Email address: wzhang@aut.ac.nz).}
\;
Je-Chiang Tsai\footnote{Department of Mathematics, National Chung Cheng University,
168, University Road, Min-Hsiung, Chia-Yi 621, Taiwan.
(Email address: tsaijc.math@gmail.com).
To whom the correspondence should be addressed (J.-C. Tsai).}
\;  and James Sneyd\footnote{Department of
Mathematics, University of Auckland, Private Bag 92019, Auckland, New Zealand
(Email address:  j.sneyd@auckland.ac.nz).}
}

\title{Curvature dependence of propagating velocity for a simplified calcium model}

\date{\today}

\maketitle
\begin{abstract}
{\small
It is known that curvature relation plays a key role in the propagation of two-dimensional waves in an excitable model.
Such a relation is believed to obey the eikonal equation for typical excitable models (e.g., the FitzHugh-Nagumo (FHN) model),
which states that the relation between the normal velocity and the local curvature is approximately linear.
In this paper, we show that for a simplified model of intracellular calcium dynamics,
although its temporal dynamics can be investigated by analogy with the FHN model,
the curvature relation does not obey the eikonal equation.
Further,
the inconsistency with the eikonal equation for the calcium model
is because of the dispersion relation
between wave speed $s$ and volume-ratio parameter $\gamma$ in the closed-cell version of the model,
not because of the separation of the fast and the slow variables as in the FHN model.
Hence this simplified calcium model may be an unexpected excitable system,
whose wave propagation properties cannot be always understood by analogy with the FHN model.
}

\bigskip

{{\bf Key Words:} calcium dynamics, eikonal equation, traveling waves, stability, FitzHugh-Nagumo model}

\bigskip
{{\bf AMS subject classifications.} 34A34, 34A12, 35K57 }

\end{abstract}

%-----------------------------------------------------------------------------
\bibliographystyle{plain}
\bibliography{autosam}

%%%%%%%%%%%%%%%%%%

\section{Introduction}
\setcounter{equation}{0}

Since the celebrated works of of Hodgkin and Huxley~\cite{Hodgkin52},
FitzHugh~\cite{Fitzhugh60} and Nagumo~\cite{Nagumo62},
wave propagation in excitable systems has been the subject of a vast number of applied mathematics studies.
In particular,
the understanding of waves in the FitzHugh-Nagumo (FHN) system has provided a great insight into waves propagation
in a wide array of biological and chemical systems~\cite{Fife79,Keener80,Tyson88,Sneyd95},
ranging from action potentials in neurons to chemical waves in the Belousov-Zhabotinskii reaction.

Although the FHN system is important in wave propagation theory of excitable systems,
not all waves in biological systems can be well understood by this classic model.
For example, earlier works on Goldbeter's model \cite{Girard92,Sneyd93,Sneyd93b},
which is derived from past theories on the mechanism underlying calcium waves and oscillations,
has revealed that Goldbeter's model is an excitable system, but the description of waves in this model
is different from that in the FHN system.
We remark that
most analytical works on Goldbeter's model are only for the reduced system which is a piecewise linear approximation to Goldbeter's model,
and which is a phenomenological formulation, and so, does not involve many biological details.

Recently, the analysis on a calcium model (CKKONS model~\cite{Champneys07,Tsai12}),
which is based on current theory for waves of intracellular calcium concentration,
indicates that the structure of one-dimensional waves in this model is quite different from that of the FHN system.
Further, the stability analysis of waves in the CKKONS model is more subtle than those for the FHN system (see \cite{Tsai12,Ghazaryan13}).
Motivated by these previous works,
one may expect that the theory of two-dimensional waves in the CKKONS model
is different from that in the FHN system.

Previous studies~\cite{Keener80,Keener86,Zykov87,Tyson88} have demonstrated that the curvature relation of waves
is crucial for the evolution of waves in two spatial dimensions.
Specifically, previous theories~\cite{Keener80,Keener86,Zykov87,Tyson88} on the classical excitable models including the FHN system
suggest that the propagation of two-dimensional waves approximately obeys the so-called eikonal equation,
which states that the relation between the normal velocity  and the local curvature is linear.
It is the curvature relation in the CKKONS model which we would like to address in this paper.
It turns out that our analysis shows that
although the temporal dynamics of the CKKONS model behaves analogously to the FHN system,
the curvature relation in the CKKONS model does not follow the eikonal equation,
and hence, the spatio-temporal behaviour of the CKKONS model cannot be understood by analogy with the FHN system.

\begin{figure}[!hbt]
\begin{center}
\vspace{0pt}
\includegraphics[width=190pt]{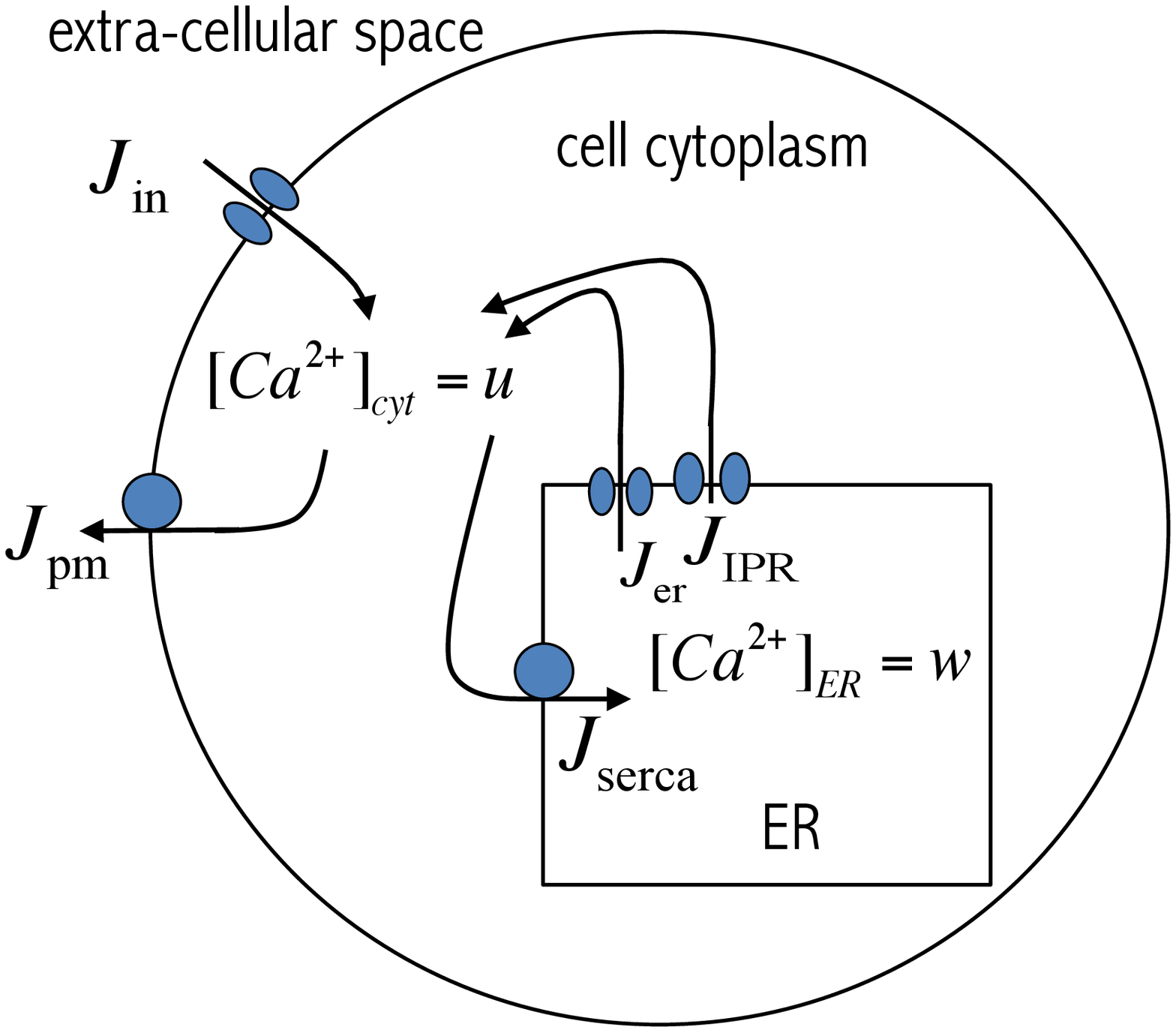}
\includegraphics[width=190pt]{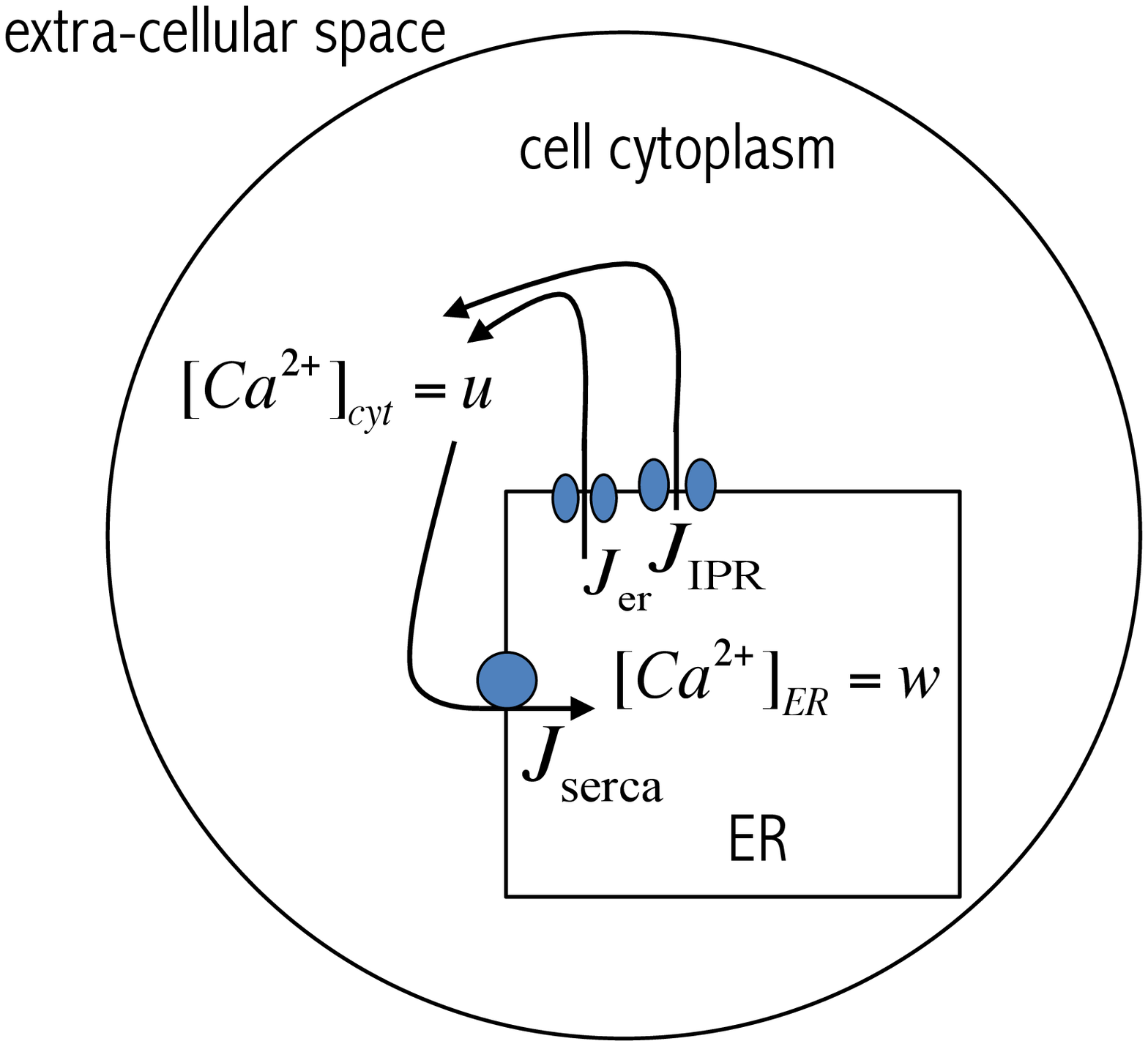}
\end{center}
\caption{\small{Schematic diagrams of the open-cell (left) and closed-cell (right) models.}}
\label{figtemp01}
\end{figure}

%--------------------------------------------------------------
\subsection{Formulation of the model}

In panel (a) of Fig.~\ref{figtemp01},
we summarize typical calcium fluxes involved in the control of cytoplasmic calcium (${\rm Ca}^{2+}$) concentration.
Hence, if we use the variables $u$ and $w$ to denote the nondimensional concentration of free cytoplasmic calcium and
the nondimensional concentration of free calcium in the ER, respectively,
then the dynamics of $u$ and $w$ are described by the following equations:
\be\label{e:1.0}
\begin{split}
     \frac{\p u}{\p t}&=D\frac{\p^2 u}{\p x^2}+J_{\rm IPR}+J_{\rm er}-J_{\rm serca}
                          +\veps\big(J_{\rm in}-J_{\rm pm}\big),\\
     \frac{\p w}{\p t}&=-\gamma\big(J_{\rm IPR}+J_{\rm er}-J_{\rm serca}\big),
\end{split}
\ee
where
the diffusion coefficient of $u$ is given by $D$,
and the parameter $\gamma$ is the ratio of the cytoplasmic volume to the ER volume,
and the constant $\veps$ controls the magnitude of the fluxes
across the membrane relative to the fluxes across the ER.
Further, we have
\begin{enumerate}
\item[$\bullet$]
$J_{\rm IPR}=k_f \frac{u^{2}}{u^{2}+\varphi_{1}^{2}} \cdot\frac{\varphi_{2}}{u^{}+\varphi_{2}}(w-u)$:\quad calcium flux through the ${\rm IP}_3$ receptor/calcium channel (${\rm IP_3R}$),
\item[$\bullet$]
$J_{\rm er}=\alpha(w-u)$:\quad calcium leak from the endoplasmic reticulum (ER),
\item[$\bullet$]
$J_{\rm serca}=k_s u$:\quad calcium flux through the ATPase calcium pumps on the membrane of the ER,
\item[$\bullet$]
$J_{\rm in}$:\quad the constant influx current of calcium from outside
the cell which is used as the main bifurcation parameter,
\item[$\bullet$]
$J_{\rm pm}= u$:\quad calcium flux through the plasma membrane ATPase calcium pumps.
\end{enumerate}
The typical values of the dimensionless parameter in the
model are as given in the following table.
\begin{table}[!hbt]
\begin{center}
\begin{tabular}{c|c|c|c|c|c|c|c}

   $D$ & $\alpha$ & $k_{s}$ & $k_f$  & $\varphi_1$ & $\varphi_2$ & $\gamma$ & $\veps $\\
\hline
  $0.025$ & $0.5$ &  $200$  & $200$ & $1$ & $2.0$  & $5.0$ & 0.001
\end{tabular}
\end{center}
\caption{Values of parameters for the dimensionless version of system~(\ref{eqn:calciumpdefull}).}
\label{tab:parameters1}
\end{table}

%%%%%%%%%%%%%%
For more detail of the model, we refer the readers to \cite{Champneys07,Tsai12}.
For ease of mathematical analysis, we recast system~(\ref{e:1.0}) in the following form:
\be\label{eqn:calciumpdefull}
\begin{split}
     \frac{\p u}{\p t}&=D\frac{\p^2 u}{\p x^2}+F(u,w) +\veps(J_{\rm in}-u),\\
     \frac{\p w}{\p t}&=-\gamma\cdot F(u,w).
\end{split}
\ee
with
\beaa
     F(u,w)&=&f(u)(w-u)-k_su:=f(u)w-g(u),\\
     f(u)&=&\alpha+k_{f}\frac{u^{2}}{u^{2}+\varphi_{1}^{2}} \cdot\frac{\varphi_{2}}{u^{}+\varphi_{2}},\\
     g(u)&=&f(u)u+k_s u.
\eeaa

\be\label{eqn:calciumpdefull-2}
\begin{split}
     \veps\frac{\p u}{\p t}&=\veps^2\frac{\p^2 u}{\p x^2}+F(u,w) +\veps(J_{\rm in}-u),\\
     \veps\frac{\p w}{\p t}&=-\gamma\cdot F(u,w).
\end{split}
\ee

%%%%%%%%%%%%%%%
%--------------------------------------------------------------
\subsection{The closed-cell model}

For many studies of calcium models,
the dynamics of the open-cell case crucially depends on that of the closed-cell case.
As we will see in this paper, there is no exception here.
For this, we state the closed-cell version of system~(\ref{eqn:calciumpdefull}).
To do this, we ignore the term for the flux of calcium across the plasma membrane in system~(\ref{eqn:calciumpdefull}).
Then the model for the closed-cell case reads (c.f. panel (b) of Fig.~\ref{figtemp01}):
\be\label{eqn:calciumpdeclosed}
\begin{split}
    \frac{\p u}{\p t}&=D\frac{\p^2 u}{\p \xi^2} + F(u,w)\\
     \frac{\p w}{\p t}&=-\gamma\cdot F(u,w).
\end{split}
\end{equation}

\medskip

Finally, the outline of the paper is as follows:
In Sec.~2, we study the temporal behavior of the open-cell model~(\ref{eqn:calciumpdefull}).
Sec.~3 is devoted to the investigation of the curvature relation of the open-cell model~(\ref{eqn:calciumpdefull}).
In particular, we use the dependence of waves of the closed-cell mode~(\ref{eqn:calciumpdeclosed})
on the volume-ratio parameter $\gamma$
to give a detailed description of the curvature relation of the closed-cell mode~(\ref{eqn:calciumpdeclosed}).
Then numerical computations indicate that for small $\veps$, the curvature relation of the open-cell model~(\ref{eqn:calciumpdefull})
is well approximated by part, not all, of the curvature relation of the closed-cell model~(\ref{eqn:calciumpdeclosed}).
In Sec.~4, we give a brief conclusion and discussion.
Finally, we reproduce the curvature relation of the FHN system in the appendix.

%%%%%%%%%%%%%%%%%%%%%%%%%%%%%%%%%%%%%%%%%%%%%%%%%%%%%%%%%%%%%%%%%%%%%%%%%%%%%%%%%%%%%%%%%%%%%%%%%%%%%%%%%%%%%%%

\section{Temporal behavior of the model}\label{sec:assumptions}
\setcounter{equation}{0}

In this section, we will show that the kinetics of the model (i.e., system~\eqref{eqn:calciumpdefull} without diffusion)
exhibits excitable behavior, and, with the introduction of the new variable $q$, its dynamics can be understood by that of the
classical excitable system--the FHN system.

To begin with, we study the equilibria of \eqref{eqn:calciumpdefull}
in the case there is no diffusion,
i.e., the equilibria of
\be\label{e:restpt-etemp20}
\begin{split}
     \frac{du}{dt}&=F(u,w)+\veps(J_{in} - u),\cr
     \frac{dw}{dt}&=-\gamma\cdot F(u,w).
     \end{split}
\ee

%-----------------------------------------------------------------------------
\subsection{Equilibria analysis: the closed-cell case}

We first deal with the closed-cell case for system~(\ref{e:restpt-etemp20}).
Hence we consider the following model:
\be\label{e:restpt-etemp25}
\begin{split}
     \frac{du}{dt}&=F(u,w),\cr
     \frac{dw}{dt}&=-\gamma\cdot F(u,w).
     \end{split}
\ee
We remark that the results in this subsection are already derived in our previous work~\cite{Tsai12}.
Since we will use the assumptions and the notations in \cite{Tsai12} for our analysis,
we restate it here for the convenience of the readers.

The equilibria of system~(\ref{e:restpt-etemp25}) lie on the curve $\Gamma_N$ in the $(u,w)$ phase space defined by
$$
  \Gamma_N:\; F(u,w)=f(u)w-g(u)=0,
$$
which can be rewritten as
$$
w=H(u):=\frac{g(u)}{f(u)},
$$
and which is the set of intersection points of the following two curves:
\beaa
      \Gamma_{1,w}&:& y=f(u)w,\\
      \Gamma_2&:& y=g(u).
\eeaa
The curve $\Gamma_N$ is
simultaneously the $u$-nullcline and the $w$-nullcline for system~(\ref{e:restpt-etemp20}).
Now we make the following assumption to ensure that $\Gamma_N$ is N-shaped.

%%%%%%%%%%%%%%%%
\begin{figure}[!hbt]
\begin{center}
\vspace{0pt}
%\resizebox{6.4in}{!}
{\includegraphics[width=410pt]{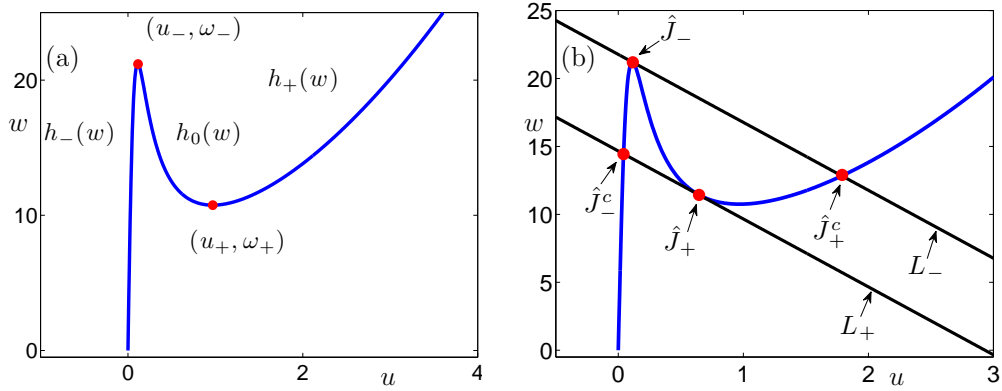}}
\end{center}
\caption{
Panel (a) shows the graph of of the nullcline $\Gamma_N$ defined by $F(u,w)=0$ or $w=H(u)$ in the $(u,w)$-plane.
Panel (b) gives the graphs of the lines $L_\pm$.
The labels on the points and curve sections in both panels are explained in {\bf Assumption (A1)-(A2)}.
The parameter values used in this figure, and in all the figures in this paper, are
as in Table~\ref{tab:parameters1}.
}
\label{nullcline}
\end{figure}
%%%%%%%%%%%%%%%%%

\medskip
\noindent
{\bf Assumption (A1):}
There exist two positive numbers $\omega_\pm$ such that
\begin{enumerate}
\item[(i)]
$\Gamma_{1,w}$ intersects $\Gamma_2$ at one point $(h_-(w),w)$ for $w\in(0,\omega_+)$;
%(see panel (a) of Fig.~\ref{intercept1});
\item[(ii)]
$\Gamma_{1,w}$ intersects $\Gamma_2$ at three points $(h_-(w),w)$, $(h_0(w),w)$ and $(h_+(w),w)$ for $w\in(\omega_+,\omega_-)$;
%(see panel (b) of Fig.~\ref{intercept1});
\item[(iii)]
$\Gamma_{1,w}$ intersects $\Gamma_2$ at one point $(h_+(w),w)$ for $w\in(\omega_-,\infty)$;
%(see panel (c) of Fig.~\ref{intercept1}).
\end{enumerate}
\medskip

We remark that for the specific choice of parameters in Table~\ref{tab:parameters1},
the model~(\ref{e:restpt-etemp25}) satisfies assumption {\bf(A1)}.
Under assumption ${\bf(A1)}$, the nullcline $\Gamma_N$ is divided into three parts:
$u=h_-(w)$, $u=h_0(w)$, and $u=h_+(w)$; we denote by $(u_-,\omega_-)$
(resp., $(u_+,\omega_+)$)  the intersection point of
the curves $u=h_-(w)$ and $u=h_0(w)$
(resp., $u=h_+(w)$ and $u=h_0(w)$).
Further, we can verify \cite{Tsai12} that
$$
\frac{dh_\pm}{dw}>0\;\mbox{ for }\;  w\in(0,\omega_+)\cup(\omega_-,\infty),
\;\mbox{and }\;
\frac{dh_0}{dw}<0 \;\mbox{ for }\;  w\in(\omega_+,\omega_-).
$$
Now we parametrise the nullcline $\Gamma_N$ with the parameter $J$, i.e., for each $J>0$,
denote by $\hE^J$ the equilibrium
$$
  \hE^J=(J,w^J):=(J,H(J)).
$$
The Jacobian matrix for system~(\ref{e:restpt-etemp20}), evaluated at $\hE^J$, is
$$
\mathcal{E}^J=\left[
\begin{array}{cc}
f'(J)w^J-g'(J) & f(J) \cr
-\gamma\big(f'(J)w^J-g'(J)\big) & -\gamma f(J)
\end{array}
\right]
$$
whose eigenvalues are $0$ and
\[\label{e:restpt-etemp27}
R(J):=f'(J)w^J-g'(J)-\gamma f(J).
\]
That  $\mathcal{E}^J$ has a zero eigenvalue is due to the fact that there is
a whole curve of equilibria;
the eigenvector corresponding to the zero eigenvalue
is the tangent vector of the curve $F(u,w)=0$ at the point $\hE^J$.
It is not clear to see the sign of the function $R(J)$, and hence the stability of $\hE^J$.
For this,
throughout the remainder of this paper, we impose the second assumption as follows.

\medskip
\noindent
{\bf Assumption (A2):}
There exist two positive numbers, $\hJ_- = \hJ_-(\gamma)$ and $\hJ_+ = \hJ_+(\gamma)$ with  $\hJ_-<\hJ_+$, such that
$R(J)<0$ for $J\in(0,\hJ_-)\cup(\hJ_+,\infty)$, $R(J)>0$ for $J\in(\hJ_-,\hJ_+)$,
and $R(\hJ_\pm)=0$ (c.f. panel (a) of Fig.~\ref{Hopf-curve}).
\medskip

We remark that for the specific choice of parameters in Table~\ref{tab:parameters1},
the model~(\ref{e:restpt-etemp25}) satisfies assumption {\bf(A1)}.
Then we have the following proposition.
\begin{proposition}\label{p:restpt-ptemp10}
  Under assumptions {\bf(A1)} and {\bf(A2)},
$\hE^J$ is a stable point of system~(\ref{e:restpt-etemp25}) for $J\in(0,\hJ_-)\cup(\hJ_+,\infty)$,
and $\hE^J$ is an unstable point of system~(\ref{e:restpt-etemp25}) for $J\in(\hJ_-,\hJ_+)$.
\end{proposition}

Now fix a $J_l\in(0,u_-)$ be such that $H(J_l) > \omega_+$.
Let $L_\gamma$ be the line in the $(u,w)$ plane defined by $u+\tfrac{1}{\gamma}w=J + \tfrac{1}{\gamma}H(J_l)$,
and define
\begin{equation}\label{e:gamma_M}
\gamma_M = \gamma_M(J_l) :=\sup\big\{\gamma > 0 |\; \mbox{$L_\gamma$ intersects $w=H(u)$ at three points} \big\}.
\end{equation}
Finally,
let $L_{\gamma,\pm}$ be the lines in the $(u,w)$ plane defined by $u+\tfrac{1}{\gamma}w=\hJ_\pm + \tfrac{1}{\gamma}H(\hJ_\pm)$,
and let $\hJ_\mp^c = \hJ_\mp^c(\gamma)$ be the $u$-coordinates of the intersection points of $L_{\gamma,\pm}$ and $u=h_\mp(w)$,
as shown in the right panel of Fig.~\ref{nullcline}.
Then we have the following lemma.
\begin{lemma}\label{l:restpt-ltemp10}
  Let $J_l\in(0,u_-)$ be such that $H(J_l) > \omega_+$,
and  $L$ be the line in the $(u,w)$ plane defined by $u+\tfrac{1}{\gamma}w=J_l + \tfrac{1}{\gamma}H(J_l)$ with $\gamma\in(0,\gamma_M)$.
Then $L$ intersects $w=H(u)$ at exactly three points:
$(J_l,H(J_l))$,
$(J_m,H(J_m))$ with $J_m\in(\hJ_-,\hJ_+)$,
and $(J_r,H(J_r))$ with $J_r\in(\hJ_+,\hJ_+^c)$.
\end{lemma}

%-----------------------------------------------------------------------------
\subsection{Equilibria analysis and temporal dynamics for the open-cell case}
\subsubsection{Equilibria analysis}

Now we turn to determine the equilibria of the open-cell model~(\ref{e:restpt-etemp20}).
The equilibrium solution lies in the $(u,w)$-plane at the intersection of the curve $J_{in} - u=0$ and the curve $\Gamma_N$.
Rearranging the functions of the curves, we then obtain the explicit  expressions of $u$ and $w$ as:
$$
u=J:=J_{in}, \quad w=H(u):=\frac{g(u)}{f(u)}.
$$
Hence for each choice of $J_{in}$ and so, for a particular choice of $J$,
there is a unique equilibrium point of system~(\ref{e:restpt-etemp20})
which is exactly $\hE^J$, and which is parametrized by the parameter $J$.

\begin{figure}[!hbt]
\begin{center}
\vspace{0pt}
{\includegraphics[width=410pt]{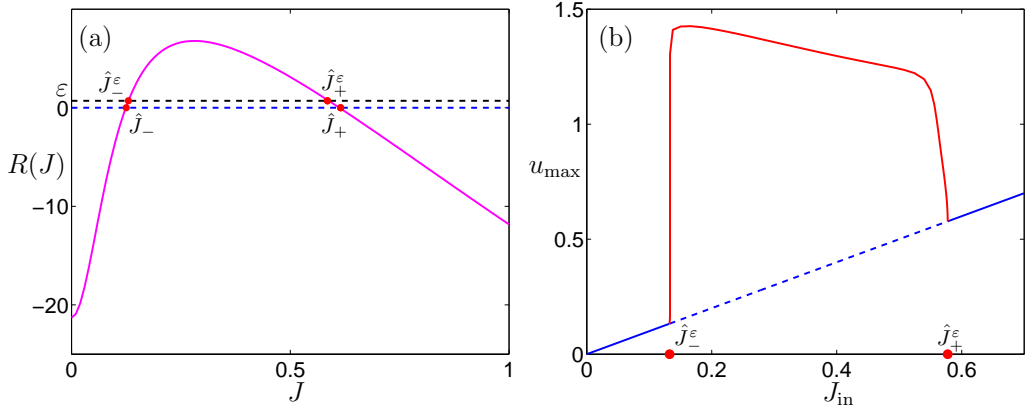}}
\end{center}
\caption{
A plot of the function $R(J)$ (panel (a)), and
and the branch of periodic orbits connecting the Hopf bifurcation points (panel (b)).
}
\label{Hopf-curve}
\end{figure}

The Jacobian matrix for system~(\ref{e:restpt-etemp20}), evaluated at $\hE^J$, is
$$
\mathcal{E}^{J,\veps}=\left[
\begin{array}{cc}
f'(J)w^J-g'(J)-\veps  & f(J) \cr
-\gamma\big(f'(J)w^J-g'(J)\big) & -\gamma f(J)
\end{array}
\right].
$$
The determinant and the trace of the Jacobian matrix are
$$
\veps \gamma f(J)
\mbox{ and }
R^{\veps}(J):=R(J) -\veps =f'(J)w^J-g'(J)-\gamma f(J)-\veps,
$$
respectively (c.f. panel (a) of Fig.~\ref{Hopf-curve}).
Hence the stability of the equilibrium point $\hE^J$ is determined by the roots of the characteristic equation
$$
    \lambda^2 - R^{\veps}(J) \lambda   + \veps \gamma f(J) = 0.
$$
According to assumption {\bf(A2)} and the relation $R^{\veps}(J) = R(J) - \veps$,
we can conclude that for each small $\veps>0$,
there exist two positive numbers, $\hJ_-^\veps$ and $\hJ_+^\veps$ with  $\hJ_- < \hJ_-^\veps < \hJ_+^\veps <\hJ_+$, such that
$R^\veps(J)<0$ for $J\in(0,\hJ_-^\veps)\cup(\hJ_+^\veps,\infty)$, $R^\veps(J)>0$ for $J\in(\hJ_-^\veps,\hJ_+^\veps)$,
and $R(\hJ_\pm^\veps)=0$. Moreover, $\hJ_\pm^\veps \to \hJ_\pm$ as $\veps\to 0$.

Together with the fact that the term $\veps \gamma f(J)$ is always greater than zero for $J>0$,
we can deduce that for $J\in(0,\hJ_-^\veps)\cup(\hJ_+^\veps,\infty)$, $\hE^J$ is a stable point of system~(\ref{e:restpt-etemp20}),
while for $J\in(\hJ_-^\veps,\hJ_+^\veps)$, $\hE^J$ is an unstable point of system~(\ref{e:restpt-etemp20}).
Further, there is a Hopf bifurcation at the critical parameters $J=\hJ_\pm^\veps$.

%%%%%%%%%%%%%%%%
%%%%%%%%%%%%%%%%%

%------------------------------------------------------------------------------------------------------------
\subsubsection{Temporal dynamics}

\begin{figure}[!hbt]
\begin{center}
\vspace{0pt}
\includegraphics[width=410pt]{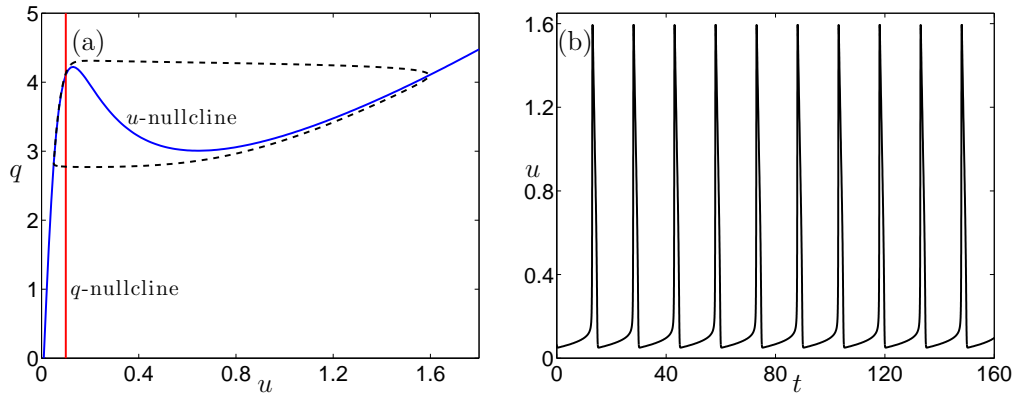}
\end{center}
\caption{\small{
(a). The phase portrait of a periodic solution of system~(\ref{e:restpt-etemp20}) is plotted in a dashed curve in the $(u,q)$-plane.
The $u$-nullcline and the $q$-nullcline are plotted in blue and red solid lines, respectively.
(b). The periodic solution in panel (a) is plotted as a function of time $t$.
The solution has fast spikes separated by longer latent periods.
}}
\label{figtemp03}
\end{figure}

Both of the Hopf bifurcations are supercritical and a stable periodic orbit is produced at the Hopf bifurcations as the parameter $J$ is varied. The numerical results of the bifurcation diagram is shown in panel (b) of Fig.~\ref{Hopf-curve}. The typical phase portrait of the periodic orbits created in the Hopf bifurcation is shown in panel (a) of Fig.~\ref{figtemp03}.
The same periodic solution is shown in the time domain in panel (b) of Fig.~\ref{figtemp03}, which shows a typical relaxation oscillation with fast spikes separated by longer latent periods.

\medskip
{\bf Insert: Biological interpretation of relaxation oscillation ???}
\medskip

The existence of a relaxation oscillation indicates that system~(\ref{e:restpt-etemp20}) is excitable.
Indeed, for $J\in(0,\hJ_-^\veps)\cup(\hJ_+^\veps,\infty)$,
a further numerical attempt suggests that superthreshold perturbations from the equilibrium point $\hE^J$ will experience a large transient excursion before returning to $\hE^J$,
while subthreshold perturbations from $\hE^J$ will decay to $\hE^J$ exponentially fast.
Hence system~(\ref{e:restpt-etemp20}) is excitable.
However, unlike the typical excitable system,
one cannot identify the fast and slow variables of system~(\ref{e:restpt-etemp20}).

To solve this puzzle, we introduce the new variable $q:=u+\frac{w}{\gamma}$.
In view of the definition of the parameter $\gamma$,
the variable $q$ can be interpreted as the total calcium concentration.
Then with a straightforward computation,
system~(\ref{e:restpt-etemp20}) can be transformed into the following system:
\be\label{eqn:calciumsin}
\begin{split}
     \frac{\di u}{\di t}&=\chi(u,q) +\veps(J_{\rm in}- u),\\
     \frac{\di q}{\di t}&=\veps(J_{\rm in}- u),
\end{split}
\ee
with
$$
     \chi(u,q):=\gamma(q-u)\cdot f(u)-g(u).
$$
In the setting of system~(\ref{eqn:calciumsin}),
it is clear that $u$ is the fast variable, while $q$ is the slow variable.
Moreover, the $u$-nullcline is $N$-shaped and the $q$-nullcline is vertical in the ($u$,$q$) phase plane.
Hence the configuration of the null clines of system~(\ref{eqn:calciumsin})
is similar to that for the FHN system.
Thus, we can conclude that the dynamics of system~(\ref{eqn:calciumsin}), and hence that of system~(\ref{e:restpt-etemp20}),
can be understood by the theory of the FHN system.

%---------------------------------------------------------------------------------------

\section{Curvature relation of waves}
\setcounter{equation}{0}

In the previous section,
we see that with the introduction of the variable $q$ which is the total calcium concentration,
the dynamics of system~\eqref{eqn:calciumpdefull} without diffusion
can be understood by the theory of the FHN system.
However, when diffusion is present,
for system~\eqref{eqn:calciumpdefull}, the separation of time scales between the variables are not so clear,
which suggests that the theory of the FHN system cannot apply to the study of waves
in system~\eqref{eqn:calciumpdefull}, at least directly.
In this section, we will analyze the curvature relation of waves in system~\eqref{eqn:calciumpdefull}
which no longer obeys the eikonal equation as $\veps\to 0$.
Since the dispersion relation and the curvature relation can determine two-dimensional waves,
and two-dimensional waves of the FHN system propagate according to the eikonal equation as $\veps\to 0$,
we may conclude that system~\eqref{eqn:calciumpdefull} is an excitable system,
but is essentially different from the FHN system.

%-------------------------------
\subsection{One-dimensional waves of the closed-cell model}\label{sec:3.1}

To begin with, we analyze the closed-cell case.
For the convenience of the readers, we restate the closed-cell model~(\ref{eqn:calciumpdeclosed}) here.
\be\label{eqn:calciumtemp30}
\begin{split}
    \frac{\p u}{\p t}&=D\frac{\p^2 u}{\p \xi^2} + F(u,w)\\
     \frac{\p w}{\p t}&=-\gamma\cdot F(u,w).
\end{split}
\end{equation}
Since the computation of the curvature relation for system~(\ref{eqn:calciumtemp30}) is based on one-dimensional waves,
we first analyze the traveling waves of system~(\ref{eqn:calciumtemp30}).
To seek a traveling wave of system~(\ref{eqn:calciumtemp30}),
we define the moving coordinate, $\xi=x+st$, where $s = s(\gamma)$ is the wave speed.
Then a traveling wave solution $(u,w)$ of system~(\ref{eqn:calciumtemp30}) will be a function of $\xi$ alone,
and its governing equation reads
\be\label{eqn:calciumtemp31}
\begin{split}
    s \frac{\di u}{\di \xi}&=D\frac{\di^2 u}{\di \xi^2} + F(u,w)\\
    s \frac{\di w}{\di \xi}&=-\gamma\cdot F(u,w),
\end{split}
\end{equation}
together with the boundary conditions
\begin{subequations}\label{e:bc}
\bea
    (u,w)\to \hE^{J_l} \;\mbox{ as }\; \xi\to-\infty,\label{e:bca}\\
    (u,w)\to \hE^{J_r} \;\mbox{ as }\; \xi\to\infty,\label{e:bcb}
\eea
\end{subequations}
for some choice of $J_l$ and $J_r$.
If $\hE^{J_l}\not= \hE^{J_r}$,
then the traveling wave solution of system~(\ref{eqn:calciumtemp30}) is a {\em traveling front}
connecting the equilibrium $\hE^{J_l}$ to the equilibrium $\hE^{J_r}$,
while if $\hE^{J_l}=\hE^{J_r}$,
then the traveling wave solution of system~(\ref{eqn:calciumtemp30}) is a {\em traveling pulse}
corresponding to an orbit of system~(\ref{eqn:calciumtemp31}) homoclinic to $\hE^{J_l}$.

\medskip

%We collect some theorems regarding traveling waves of system~{\rm(\ref{eqn:calciumtemp30})} from \cite{Tsai12}.

Following the argument of \cite[Proposition 3]{Tsai12},
we can derive the following proposition on the
necessary condition for the existence of traveling waves of system~(\ref{eqn:calciumtemp30}).
\begin{proposition}\label{p:neceptemp5}
Let $J_l\in(0,u_-)$ be such that $H(J_l) > \omega_+$.
If system~{\rm(\ref{eqn:calciumtemp30})} admits a traveling front or pulse with positive wave speed,
then  $\gamma\in(0,\gamma_M)$.
\end{proposition}

Conversely, for the existence of traveling waves of system~{\rm(\ref{eqn:calciumtemp30})},
we fix a $J_l\in(0,u_-)$ such that $H(J_l) > \omega_+$, and let $\gamma\in(0,\gamma_M)$.
Then Lemma~\ref{l:restpt-ltemp10} guarantees that there exists a unique $J_r\in(\hJ_+,\hJ_+^c)$
such that the following relation holds.
$$       %\be\label{e:tfastfrontetemp10}
  {J_r}+\frac{w^{J_r}}{\gamma}=J_l+\frac{w^{J_l}}{\gamma}.
$$
Now if we follow the arguments of \cite[Theorem 1-3]{Tsai12},
we have the following theorems for the existence of traveling waves of system~(\ref{eqn:calciumtemp30}).
\begin{theorem}\label{t:wave_ex}

\begin{enumerate}
\item[\rm(a)]
System~{\rm(\ref{eqn:calciumtemp30})} admits a unique {\rm(}up to a translation{\rm)} traveling front $(u,w)$
with wave speed $s_F=s_F(\gamma)>0$ which satisfies the following properties:
\begin{enumerate}
\item[{\rm(i)}]
 $(u(-\infty),w(-\infty))=\hE^{J_l}$ and $(u(\infty),w(\infty))=\hE^{J_r}$;
\item[{\rm(ii)}]
 $u'>0$ on $\BbR$.
\end{enumerate}
\item[\rm(b)]
System~{\rm(\ref{eqn:calciumtemp30})} admits a unique {\rm(}up to a translation{\rm)} traveling front $(u,w)$
with wave speed $s_B=s_B(\gamma)>0$ which satisfies the following properties:
\begin{enumerate}
\item[{\rm(i)}]
 $(u(-\infty),w(-\infty))=\hE^{J_r}$ and $(u(\infty),w(\infty))=\hE^{J_l}$;
\item[{\rm(ii)}]
 $u' < 0$ on $\BbR$.
\end{enumerate}
\item[\rm(c)]
\begin{enumerate}
\item[{\rm(i)}]
 If $s_F(\gamma)\le s_B(\gamma)$, then there does not exist a traveling pulse $(u,w)$ of system~{\rm(\ref{eqn:calciumtemp30})}
with positive wave speed and $(u,w)(\pm\infty)=\hE^{J_l}$.
\item[{\rm(ii)}]
 If $s_F(\gamma)>s_B(\gamma)$, then system~{\rm(\ref{eqn:calciumtemp30})} admits a unique
{\rm(}up to a translation{\rm)} traveling pulse $(u,w)$ with positive wave speed $s_P(J_l)\in(s_B(J_l),s_F(J_l))$
and $(u,w)(\pm\infty)=\hE^{J_l}$.
\end{enumerate}
\end{enumerate}
\end{theorem}

%{\rm [Existence of traveling fronts with monotone increasing profiles for the $u$-component]}

{\em Throughout this paper,
we will call the traveling wave solutions established
in assertions {\rm(i)}, {\rm(ii)}, and {\rm(iii)} of Theorem~\ref{t:wave_ex}
as the wave front, the wave back, and the wave pulse, respctively.
We also retain the notations: $s_F(\gamma)$, $s_B(\gamma)$, and $s_P(\gamma)$.
}

\begin{figure}[!hbt]
\begin{center}
\vspace{0pt}
\includegraphics[width=250pt]{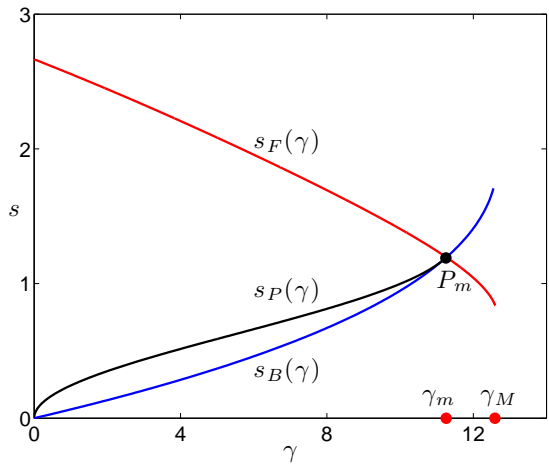}
\end{center}
\caption{Numerically computed dependence of the wave speed functions $s_F(\gamma)$ (red curve), $s_B(\gamma)$ (blue curve), and $s_P(\gamma)$
(black curve) on the parameter $\gamma\in(0,\gamma_M)$, for system~(\ref{eqn:calciumtemp30}) with parameter values as in Table~\ref{tab:parameters1}.
Throughout this paper, the left-hand side boundary condition for waves $(u,v)$ are set to be $u(-\infty)=J_l=0.06$.
Traveling pulses do not exist for $\gamma\in(0,\gamma_m)$,
and so the associated wave speed $s_P(\gamma)$ is not defined for those $\gamma$.
}
\label{figtemp2}
\end{figure}

For system~\eqref{eqn:calciumtemp31},
we fix a  $J_l\in(0,u_-)$ such that $H(J_l) > \omega_+$, and then numerically compute wave front, wave back
and wave pulse solutions in the $(\gamma,s)$-plane using the software package AUTO~\cite{auto}.
Fig.~\ref{figtemp2} plots wave speed $s$ as a function of parameter $\gamma$ for three types of traveling waves.
Numerical evidence indicates that waves cease to exist for $\gamma > \gamma_M$,
which is predicted by Proposition~\ref{p:neceptemp5}.
We denote by $P_{\rm m}:=(\gamma_{\rm m},s_F(\gamma))$ the intersection point of the curves $s_F(\gamma)$ and $s_B(\gamma)$.
For $\gamma\in(0,\gamma_{\rm m})$, the numerical computation indicates $s_F(\gamma)>s_B(\gamma)$,
and hence part (c)-(ii) of Theorem~\ref{t:wave_ex} suggests the existence of wave pulses of system~(\ref{eqn:calciumtemp31}),
while for $\gamma\in(\gamma_{\rm m},\gamma_{\rm M})$, we have $s_F(\gamma)\le s_B(\gamma)$,
and so part (c)-(i) of Theorem~\ref{t:wave_ex} prevents the existence of traveling pulses.
The numerical results for the existence of wave pulses are consistent with these predictions.

%-------------------------------
\subsection{Two-dimensional waves of the closed-cell model}\label{sec:3.2}

To determine the curvature dependence of the propagating normal velocity of two-dimensional waves
for system~(\ref{eqn:calciumtemp30}),
we follow \cite{Zykov87,Tyson88,Meron92}) to assume that compared with the radius of curvature,
the wave is very thin so that its front and rear side have the same curvature.
Under this assumption,
the wave profile $(u,w)$ of system~(\ref{eqn:calciumtemp30}) with curvature $\kappa$ and normal velocity $\ts$
satisfies the equations (see \cite{Zykov87,Tyson88,Meron92}):
\be\label{eqn:calciumtemp50}
\begin{split}
    (\ts + D\kappa) \frac{\di u}{\di \xi}&=D\frac{\di^2 u}{\di \xi^2}  + F(u,w),\\
    \ts \frac{\di w}{\di \xi}&=-\gamma_0\cdot F(u,w),
\end{split}
\ee
where $\xi = x + \ts t$ is the moving coordinate and $\gamma_0$ is fixed to be equal to 5.0 for numerical simulation.
We remark that
system~(\ref{eqn:calciumtemp50}) is the governing equations for one-dimensional traveling waves
of the following system:
\be\label{eqn:calciumtemp40}
\begin{split}
     \frac{\p u}{\p t}&=D\frac{\p^2 u}{\p x^2} - D\kappa\frac{\p u}{\p x} + F(u,w),\\
     \frac{\p w}{\p t}&=-\gamma_0\cdot F(u,w).
\end{split}
\ee
Hence the curvature dependence of normal velocity for system~(\ref{eqn:calciumtemp30})
is given by the relation between the wave speed $\ts$ of traveling waves of system~(\ref{eqn:calciumtemp40}) and the curvature parameter $\kappa$.

In below, we will borrow the idea of \cite{Zykov87,Zykov98}
to describe a method which can determine the dependence of the wave speed $\ts(\kappa)$ on $\kappa$
using the wave speed function $s(\gamma)$ associate with system~\eqref{eqn:calciumtemp30}.
To see this, by multiplying the second equation of system~(\ref{eqn:calciumtemp50}),
and setting
\be\label{eqn:calciumtemp55}
  \gamma^* = \gamma_0 \cdot\frac{(\ts + D\kappa)}{\ts}
\;\mbox{ and  }\;
   s^* = \ts + D\kappa,
\ee
system~(\ref{eqn:calciumtemp50}) can be written as the following
\be\label{eqn:calciumtemp71}
\begin{split}
    s^* \frac{\di u}{\di \xi}&=D\frac{\di^2 u}{\di \xi^2} + F(u,w),\\
    s^* \frac{\di w}{\di \xi}&=-\gamma^* F(u,w),
\end{split}
\end{equation}
which is the same as system~\eqref{eqn:calciumtemp30} with $(s,\gamma)$ replaced by $(s^*,\gamma^*)$.
Hence $s^* = s(\gamma^*)$ which is the wave speed of traveling waves of system~\eqref{eqn:calciumtemp30} with $\gamma=\gamma^*$.
Note that $s(\gamma^*)$ can be $s_F(\gamma^*)$, $s_B(\gamma^*)$, or $s_P(\gamma^*)$.
Using (\ref{eqn:calciumtemp55}) to express $s^*$ in terms of $\gamma^*$,
we have
\be\label{eqn:calciumtemp77}
   s^* = D\kappa \frac{\gamma^*}{\gamma^*-\gamma_0}.
\ee
Together with the relation $s^* = s(\gamma^*)$, we arrive at the equation
\be\label{eqn:calciumtemp90}
   s(\gamma^*)=D\kappa \frac{\gamma^*}{\gamma^*-\gamma_0}.
\ee
For each $\kappa$,
the solution of \eqref{eqn:calciumtemp90} can be understood by simultaneously plotting the function $s(\gamma^*)$ and the hyperbola function $\Phi^\kappa(\gamma^*) := D\kappa \frac{\gamma^*}{\gamma^*-\gamma_0}$ against $\gamma^*$ in the $(\gamma^*,s)$-plane.
The $\gamma^*$ value of the intersection point of two graphs of the functions $s(\gamma^*)$ and $\Phi^\kappa(\gamma^*)$ is the root $\gamma^*$ of \eqref{eqn:calciumtemp90}.
Once the root $\gamma^*$ of Eq.~\eqref{eqn:calciumtemp90} is found,
then the wave speed $\ts = \ts(\kappa)$ associated with system~(\ref{eqn:calciumtemp50})
can be determined by the second equation of (\ref{eqn:calciumtemp55}) and Eq.~(\ref{eqn:calciumtemp77}), that is,
\be\label{eqn:calciumtemp93}
   \ts(\kappa)  =  s^* - D\kappa = D\kappa \cdot  \frac{\gamma_0}{\gamma^*-\gamma_0}.
\ee

%------------------------------------------------------------------------------------------------------------
\subsection{Curvature relation of the closed-cell model}

We apply the method discussed in Sec.~\ref{sec:3.2} to understand
the curvature dependence of the normal velocity for the closed-cell model (\ref{eqn:calciumtemp30}).

\subsubsection{The case for $s(\gamma^*) = s_F(\gamma^*)$}\label{sec:3.3.1}

\begin{figure}[!hbt]
\begin{center}
\vspace{0pt}
\includegraphics[width=410pt]{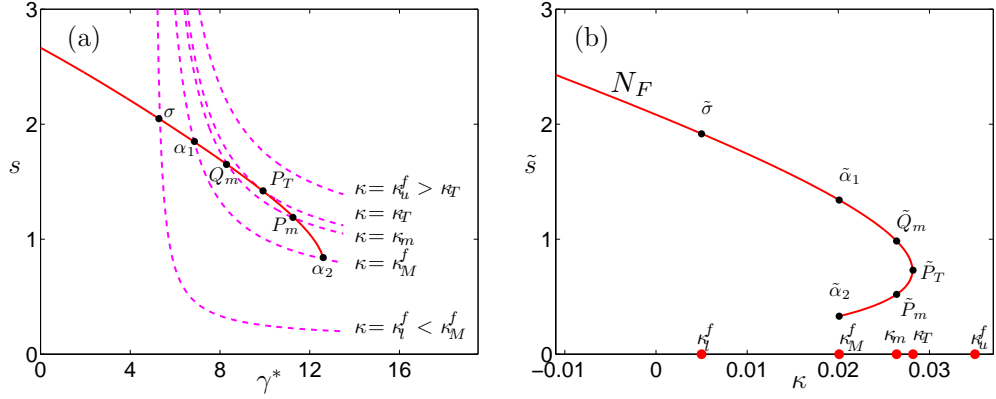}
\end{center}
\caption{\small{
(a) The function $s(\gamma^*)$ is plotted against $\gamma^*$ in the red solid curve
and the hyperbola function $\Phi^\kappa(\gamma^*)$ reprenting the right-hand side of \eqref{eqn:calciumtemp90} are plotted in the pink dashed curves for various $\kappa$ in the $(\gamma^*,s)$-plane.
(b)
The dependence of $\ts$ on $\kappa$
by numerically solving Eqs.~(\ref{eqn:calciumtemp90}) and (\ref{eqn:calciumtemp93}),
and the resulting curve in the $(\kappa,\ts)$-plane is denoted by $N_F$.
Here we choose $s(\gamma^*) = s_F(\gamma^*)$ and the parameter values from Table~\ref{tab:parameters1}.
The labels of the intersection points are explained further in the text.
}}
\label{wavefront}
\end{figure}

To begin with, we set $s(\gamma^*) = s_F(\gamma^*)$ and solve Eq.~(\ref{eqn:calciumtemp90}) geometrically
where $s_F(\gamma^*)$ is the wave speed of wave fronts of system~(\ref{eqn:calciumtemp30}) with $\gamma=\gamma^*$.
In panel (a) of Fig.~\ref{wavefront},  we plot the function $s_F(r^*)$ against $\gamma^*$ in the red solid curve, and  the hyperbola function $\Phi^\kappa(\gamma^*)$ corresponding to the right-hand side of \eqref{eqn:calciumtemp90} for $\kappa=\kappa_l^f$, $\kappa=\kappa_M^f$,
$\kappa=\kappa_m$, $\kappa=\kappa_T$, and $\kappa=\kappa_u^f$ in the pink dashed curves in the $(\gamma^*,s)$-plane.
We note that the $\gamma^*$ values associated with $\kappa=\kappa_M^f$ and $\kappa=\kappa_m$ are $\gamma_M$ and $\gamma_m$, respectively.
In panel (a) of Fig.~\ref{wavefront}, we can see that as the curvature $\kappa$ increases, the hyperbola moves up and to the right of the plot, and the number of intersection points between the graphs of the functions $s_F(\gamma^*)$ and
$\Phi^\kappa(\gamma^*)$ changes with respect to the value of $\kappa$.
Specifically, when $\kappa$ is less than the critical value $\kappa_M^f=0.0201$, say $\kappa=\kappa_l^f < \kappa_M^f$,
there is one intersection point $\sigma$ between the graphs of the functions $s_F(\gamma^*)$ and
$\Phi^\kappa(\gamma^*)$.
At the critical value $\kappa=\kappa_M^f$, the hyperbola intersects the graph of the function $s_F(\gamma^*)$ at two points $\alpha_1$ and $\alpha_2$.
One of the intersection points $\alpha_2$ is an endpoint of the graph of the function $s_F(\gamma^*)$.
When $\kappa$ is greater than $\kappa_M^f$ and less than another critical value $\kappa_T=0.0282$,
there are two intersection points between  the graphs of the functions $s_F(\gamma^*)$ and
$\Phi^\kappa(\gamma^*)$.
At $\kappa=\kappa_T$, the hyperbola is tangent to the graph of the function $s_F(\gamma^*)$ at the point $P_T$.
The graphs of the functions $s_F(\gamma^*)$ and
$\Phi^\kappa(\gamma^*)$ do not intersect for $\kappa=\kappa_u^f>\kappa_T$.

From the aforementioned discussion,
Eq.~(\ref{eqn:calciumtemp90}) admits a root if and only if $\kappa\in(-\infty,\kappa_T]$.
Then for a given curvature $\kappa\in(0,\kappa_T]$,
we can substitute the $\gamma^*$ value, which is determined by Eq.~(\ref{eqn:calciumtemp90}),
into \eqref{eqn:calciumtemp93}  to obtain the dependence $\ts$ on the curvature $\kappa$.
Applying this method to the intersection points $\sigma$, $\alpha_1$, $\alpha_2$, $Q_m$, $P_m$, and $P_T$
of the graphs of the functions $s_F(\gamma^*)$ and
$\Phi^\kappa(\gamma^*)$ in panel (a) of Fig.~\ref{wavefront}
whose corresponding $\kappa$-coordinates are $\kappa_l^f$, $\kappa_M^f$, $\kappa_M^f$, $\kappa_m$, $\kappa_m$, and $\kappa_T$, respectively,
we can locate the corresponding points $\tilde \sigma$, $\tilde \alpha_1$, $\tilde \alpha_2$, $\tilde{Q}_m$, $\tilde{P}_m$, and $\tilde{P}_T$
in the $(\kappa,\ts)$-plane, as shown in panel (b) of Fig.~\ref{wavefront}.
We note that $\tilde{P}_T$ is the turing point for the curve of curvature relation in the $(\kappa,\ts)$-plane.

From panel (b) of Fig.~\ref{wavefront},
we can conclude that
the wave speed $\ts$ is defined only for $(-\infty,\kappa_T]$, and
the dependence of the wave speed $\ts$ on $\kappa$ can be characterized as follows:
$\ts$ is single-valued for $\kappa\in(-\infty,\kappa_M^f)\cup\{k_T\}$,
and double-valued for $\kappa\in(\kappa_M^f,\kappa_T)$.
In the remaining of this paper,
we denote by $N_F$ the curve of curvature relation in the $(\kappa,\ts)$-plane
which is constructed in this subsection.

%---------------------------------------------------------------------------------
\subsubsection{The case for $s(\gamma^*) = s_B(\gamma^*)$}\label{sec:3.3.2}
\begin{figure}[!hbt]
\begin{center}
\vspace{0pt}
\includegraphics[width=410pt]{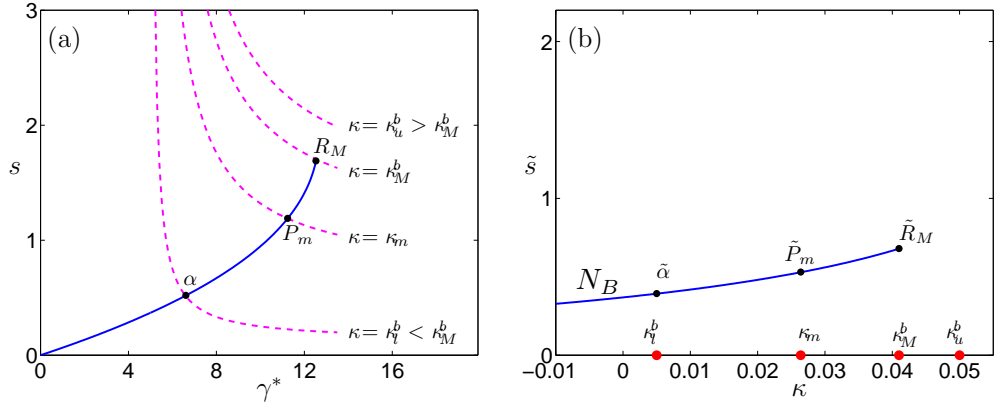}
\end{center}
\caption{\small{
(a) The function $s_B(\gamma^*)$ is plotted  in the blue solid curve
and the hyperbola function $\Phi^\kappa(\gamma^*)$  are plotted in the pink dashed curves for $\kappa=\kappa_l^b$, $\kappa=\kappa_m$, $\kappa=\kappa_M^b$, and $\kappa=\kappa_u^b$ in the $(\gamma^*,s)$-plane.
(b)
The dependence of $\ts$ on $\kappa$ for the case $s(\gamma^*) = s_B(\gamma^*)$ is found
by numerically solving Eqs.~(\ref{eqn:calciumtemp90}) and (\ref{eqn:calciumtemp93}),
and the resulting curve in the $(\kappa,\ts)$-plane is denoted by $N_B$.
The parameter values are from Table~\ref{tab:parameters1}.
The labels of the intersection points are explained further in the text.
}}
%\caption{\small{\textit{
%(a) The wave back bifurcation curve $s(\gamma^*)$ is plotted in blue solid curve and the hyperbola functions on the right hand side of \eqref{eqn:calciumtemp90} are plotted in pink dashed curves
%for $\kappa=\kappa_l^b$, $\kappa=\kappa^b$ and $\kappa=\kappa_u^b$
%in the $(\gamma^*,s)$ parameter plane. The blue solid curve and the pink dashed curves intercept at two points $\alpha$ and $\beta$.
%(b) The wave back bifurcation curve of \eqref{eqn:calciumtemp93} is plotted in blue
%in the $(\kappa,\ts)$ parameter plane. The labels of the intersection points are explained further in the text.
%}}}
\label{waveback}
\end{figure}

Similar to Sec.~3.3.1, we study the solution of Eq.~(\ref{eqn:calciumtemp90}) in the case of $s(\gamma^*) = s_B(\gamma^*)$. As defined in Sec.~3.1,
$s_B(\gamma^*)$ is the wave speed of wave backs of system~(\ref{eqn:calciumtemp30}) with $\gamma=\gamma^*$.
Panel (a) of Fig.~\ref{waveback} shows the function $s_B(r^*)$  in the blue solid curve, and  the hyperbola function $\Phi^\kappa(\gamma^*)$
%corresponding to the right-hand side of \eqref{eqn:calciumtemp90}
for $\kappa=\kappa_l^b$,
$\kappa=\kappa_m$, $\kappa=\kappa_M^b$, and $\kappa=\kappa_u^b$ in the pink dashed curves in the $(\gamma^*,s)$-plane.
Note that the $\gamma^*$ values associated with $\kappa=\kappa_M^b$ and $\kappa=\kappa_m$ are $\gamma_M$ and $\gamma_m$, respectively.
Similar to panel (a) of Fig.~\ref{wavefront}, the hyperbola moves up with  increasing $\kappa$, and the functions $s_B(\gamma^*)$ and
$\Phi^\kappa(\gamma^*)$ intersect at lower values of $\kappa$ and do not intersect at higher values of $\kappa$.
In particular, when $\kappa$ is less than the critical value $\kappa_M^b=0.0435$, say $\kappa=\kappa_l^b < \kappa_M^b$,
the graphs of the functions $s_B(\gamma^*)$ and
$\Phi^\kappa(\gamma^*)$ intersect at one point $\alpha$.
At the critical value $\kappa=\kappa_M^b$, the hyperbola touches the graph of the function $s_B(\gamma^*)$ at its end point $R_M$.
Therefore, Eq.~(\ref{eqn:calciumtemp90}) admits a zero  of $\gamma^*$ for $\kappa\in(-\infty,\kappa_M^b]$, and no root exists for  $\kappa\in(\kappa_M^b,\infty)$.
Note that the $\gamma^*$-coordinates of the points $P_m$ and $R_M$ are $\gamma_m$ and $\gamma_M$, respectively.

To summarize,
given a curvature $\kappa\in(-\infty,\kappa_M^b]$,
we can obtain the root $\gamma^*$ by solving Eq.~(\ref{eqn:calciumtemp90}), as described in the above paragraph;
then find the corresponding $\ts$ value
by inserting $\gamma^*$ into Eq.~\eqref{eqn:calciumtemp93}.
For instance, the intersection points  $\alpha$, $P_m$ and $R_m$,
marked as black dots in the $(\gamma^*,s)$-plane in panel (a) of Fig.~\ref{waveback},
%of the graphs of the functions $s_B(\gamma^*)$ and
%$\Phi^\kappa(\gamma^*)$ in panel (a) of Fig.~\ref{waveback}
%whose corresponding $\kappa$-coordinates are $\kappa_l^b$, $\kappa_m$, and $\kappa_M^b$, respectively,
are mapped to the corresponding points $\tilde \alpha$, $\tilde{P}_m$, and $\tilde{R}_m$
in the $(\kappa,\ts)$-plane, as shown in panel (b) of Fig.~\ref{waveback}.
As seen from panel (b) of Fig.~\ref{waveback},  there is one value of wave speed $\ts$ for each $\kappa \in (-\infty,\kappa_M^b]$.

%---------------------------------------------------------------------------------
\subsubsection{The case for $s(\gamma^*) = s_P(\gamma^*)$}\label{sec:3.3.3}

\begin{figure}[!hbt]
\begin{center}
\vspace{0pt}
\includegraphics[width=410pt]{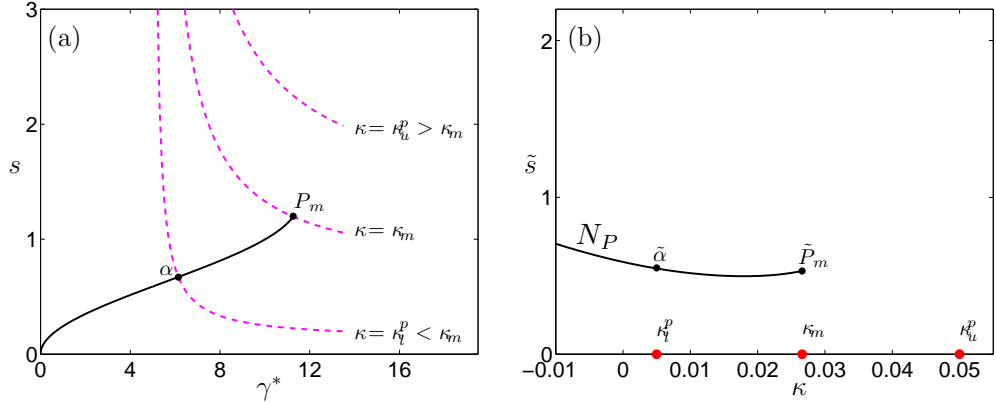}
\end{center}
\caption{\small{
(a) The function $s_P(\gamma^*)$ is plotted  in the black solid curve
and the hyperbola function $\Phi^\kappa(\gamma^*)$ on the right-hand side of \eqref{eqn:calciumtemp90} are plotted in the pink dashed curves for $\kappa=\kappa_l^p$, $\kappa=\kappa_m$ and $\kappa=\kappa_u^p$ in the $(\gamma^*,s)$-plane.
(b)
The dependence of $\ts$ on $\kappa$ for the case $s(\gamma^*) = s_P(\gamma^*)$ in the $(\kappa,\ts)$-plane,
and the resulting curve in the $(\kappa,\ts)$-plane is denoted by $N_P$.
% is found by numerically solving Eqs.~(\ref{eqn:calciumtemp90}) and (\ref{eqn:calciumtemp93}).
The parameter values are from Table~\ref{tab:parameters1}.
The labels of the intersection points are explained further in the text.
}}
%\caption{\small{\textit{
%(a) The pulse bifurcation curve $s(\gamma^*)$ is plotted in black solid curve and the hyperbola functions on the right hand side of \eqref{eqn:calciumtemp90} are plotted in pink dashed curves for various $\kappa$ in the $(\gamma^*,s)$ parameter plane.
%(b) The pulse bifurcation curve of \eqref{eqn:calciumtemp93} is plotted in black
%in the $(\kappa,\ts)$ parameter plane. The labels of the intersection points are explained further in the text.
%}}}
\label{pulse}
\end{figure}

The similar method also applies to the case for $s(\gamma^*) = s_P(\gamma^*)$,
where $s_P(\gamma^*)$ is the wave speed of wave pulses of system~(\ref{eqn:calciumtemp30}) with $\gamma=\gamma^*$.

In panel (a) of Fig.~\ref{pulse}, the function $s_P(r^*)$ is shown in the black solid curve, and the hyperbola function $\Phi^\kappa(\gamma^*)$
are shown in the pink dashed curves
for $\kappa=\kappa_p^l$,
$\kappa=\kappa_m$ and $\kappa=\kappa_p^u$  in the $(\gamma^*,s)$-plane.
Note that the $\gamma^*$ value associated with $\kappa=\kappa_m$ is $\gamma_m$.
The hyperbola moves up and to the right with  increasing $\kappa$, and the function $s_P(\gamma^*)$ and the hyperbola
$\Phi^\kappa(\gamma^*)$  intersect at lower values of $\kappa$ and do not intersect at higher values of $\kappa$.
In particular, when $\kappa$ is less than the critical value $\kappa_m=0.0266$, say $\kappa=\kappa_l^p < \kappa_m$, one intersection point, $\alpha$, exists between
the graphs of the functions $s_P(\gamma^*)$ and $\Phi^\kappa(\gamma^*)$.
At the critical value $\kappa=\kappa_m$, the hyperbola touches the graph of the function $s_P(\gamma^*)$ at its end point $P_m$.
In other words, there exists a root of $\gamma^*$ for Eq.~(\ref{eqn:calciumtemp90})  in the case $\kappa\in(-\infty,\kappa_m]$, and there exists no root  for  $\kappa\in(\kappa_m,\infty)$.

Hence, Eq.~(\ref{eqn:calciumtemp90}) has a root $\gamma^*$  for a given curvature $\kappa\in(-\infty,\kappa_m]$; and the corresponding $\ts$ value can be found by substituting the root $\gamma^*$ into Eq.~\eqref{eqn:calciumtemp93}.
Applying the idea, we can obtain the points $\tilde \alpha$ and $\tilde{P}_m$ in the $(\kappa,\ts)$-plane shown in panel (b) of Fig.~\ref{pulse} from the mapping of the intersection points  $\alpha$ and $P_m$,
marked as black dots in the $(\gamma^*,s)$-plane in panel (a) of Fig.~\ref{pulse}.
As can be seen from panel (b) of Fig.~\ref{pulse},  there is one-to-one correspondence between the wave speed $\ts$ and the curvature $\kappa$ for $\kappa \in (-\infty,\kappa_m]$.

%---------------------------------------------------------------------------------
\subsubsection{A summary and stability of waves}

\begin{figure}[!hbt]
\begin{center}
\vspace{0pt}
\includegraphics[width=250pt]{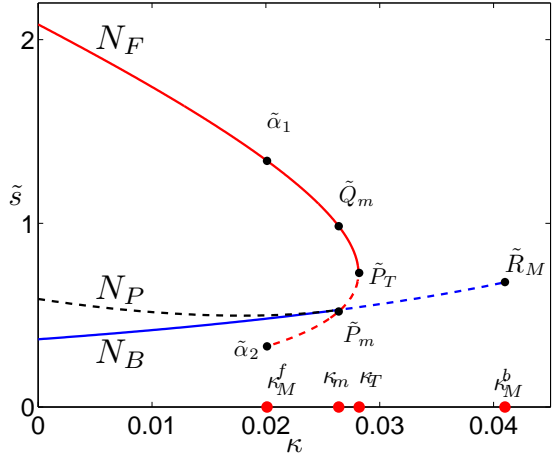}
\end{center}
\caption{\small{
The dependence of wave speed $\ts$ on curvature $\kappa$ given by the curvature relation
in the closed-cell model~(\ref{eqn:calciumtemp40}).
Such a dependence
is exactly given by the union of the curve $N_F$ (red curve) in Sec.~\ref{sec:3.3.1}, the curve $N_B$ (blue curve) in Sec.~\ref{sec:3.3.2},
and the curve $N_P$ (black curve) in Sec.~\ref{sec:3.3.3}.
Solid curves for part of $N_F$ and part of $N_B$ indicate that the corresponding waves are stable,
while dashed curves for $N_P$, part of $N_F$, and part of $N_B$ indicate that the corresponding waves are unstable.
}}
\label{stability_closed}
\end{figure}

With the results established in Sec.~\ref{sec:3.3.1}-Sec.~\ref{sec:3.3.3},
we can characterize curvature relation in the closed-cell model~\eqref{eqn:calciumtemp40}.
The results for curvature relation are depicted in Fig.~\ref{stability_closed}.
Fig.~\ref{stability_closed} plots wave speed $\ts$ against curvature parameter $\kappa$.
The curves of curvature relations in the $(\kappa,\ts)$-plane are
exactly given by the union of the curve $N_F$ (red curve) in Sec.~\ref{sec:3.3.1}, the curve $N_B$ (blue curve) in Sec.~\ref{sec:3.3.2},
and the curve $N_P$ (black curve) in Sec.~\ref{sec:3.3.3}.
We note that the $\gamma^*$ values associated with $\kappa=\kappa_M^f$, $\kappa=\kappa_m$,
$\kappa=\kappa_T$, and $\kappa=\kappa_M^b$ are $\gamma_M$, $\gamma_m$, $\gamma_T$, and $\gamma_M$, respectively.
The points in Fig.~\ref{stability_closed} are already explained in Sec.~\ref{sec:3.3.1}-Sec.~\ref{sec:3.3.3}.
The curve $N_P$ stops at the point $\tilde{P}_m$ due to the fact that
$s_F(\gamma)\le s_B(\gamma)$ for $\gamma\in[\gamma_m,\gamma_M)$ (see Sec.~\ref{sec:3.1}).

For the stability of various wave solutions along the curve of curvature relation,
we compute the spectrum of waves.
The first step of finding the waves' spectrum is to approximate system~\eqref{eqn:calciumtemp40}
by a system of time-dependent ordinary differential equations (ODEs).
We discretize the spatial variable $\xi$ to obtain $\{\xi_i\}^{n+1}_{i=1}$, where $n+1$ is the number of grid points, and define the
dependent variables $u_i (t) = u(\xi_i, t)$, and $w_i (t) = w(\xi_i, t)$,  for $i = 1, 2, 3, . . ., n, n+1$.
Periodic boundary conditions $u_1=u_{n+1}$ and $w_1=w_{n+1}$ are employed
for the computation of the spectrum~\cite{Sandstede00}.
For more information regarding the choice of the boundary conditions for the stability analysis of waves,
see \cite{Sandstede00,Beyn90}.
The distances between the grid points are all equal to $h$.

Specifically, system~\eqref{eqn:calciumtemp40} is approximated by a system of ODEs
\be\label{e:wave_stability_temp10}
\begin{split}
\frac{du_i}{dt}&=D\frac{u_{i+1}+u_{i-1}-2u_{i}}{h^2}-(\ts+D\kappa)\frac{u_{i+1}-u_i}{h}+F(u_i,w_i) \\
\frac{dw_i}{dt}&=-\ts\frac{w_{i+1}-w_i}{h}-\gamma_0 F(u_i,w_i),
\end{split}
\ee
for $i=1 ,2 ,3, ..., n$,
which, for the ease of notations, is rewritten in the form:
$$
\frac{dX}{dt} = \mathcal{G}(X),
$$
with $X = (u_1, u_2, . . ., u_n, w_1, w_2, . . ., w_n)^T \in R^{2n}$.

We can obtain a discretization $\tilde{X}$ of the traveling wave solutions of system~(\ref{e:wave_stability_temp10}) by solving $ \mathcal{G}(\tilde{X})=0$.
A good initial estimate of the solution $\tilde{X}$ can be found numerically by boundary-value solvers such as AUTO \cite{auto};
then the initial estimate can be improved to required accuracy by the Newton-Raphson iterations.
Note that we need to compute the Jacobian matrix $D_X \mathcal{G}(\tilde{X})$ in the process of the iteration step of the Newton-Raphson method.
After the solution converges, the Jacobian matrix $D_X \mathcal{G}(\tilde{X})$ is then used to obtain the spectrum of $\tilde{X}$ by solving the eigenvalue problem:
$$D_X \mathcal{G}(\tilde{X})Y=\lambda Y, \quad Y\in R^{2n}.$$
We can determine the stability of the discretization waves $\tilde{X}$ by checking whether the real part of the largest eigenvalue of the spectrum is greater than zero. If the real part of the largest eigenvalue of the spectrum is greater than zero, then the wave solution is unstable, otherwise the wave solution is stable.

     We apply the above method to determine the stability of wave solutions of system~\eqref{eqn:calciumtemp40}.
The result is summarized in Fig.~\ref{stability_closed}.

\begin{figure}[!hbt]
\begin{center}
\vspace{0pt}
\includegraphics[width=250pt]{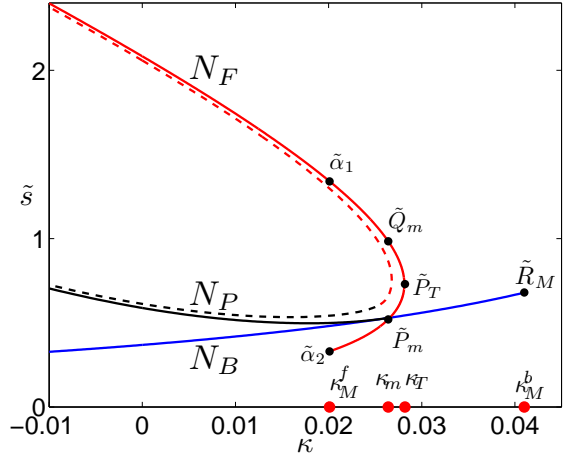}
\end{center}
\caption{\small{
Curvature relations for system~(\ref{eqn:calciumpdeclosed}) (the closed-cell model) with $J_l=0.06$
are plotted in the red and black solid curves which are part of those plotted in Fig.~\ref{stability_closed},
while
curvature relations for system~(\ref{eqn:calciumpdefull}) (the open-cell model) with $\veps=0.001$ and $J_{\rm in}=0.06$
are plotted in the red and black dashed curves.
}}\label{fig:curvature_relation_open}
\end{figure}

%------------------------------------------------------------------------------------------------------------
\subsection{Curvature relation for the the open-cell model}

In this subsection, we turn to investigate the curvature relation of the open-cell model~(\ref{eqn:calciumpdefull}).
Following the argument for the closed-cell model~(\ref{eqn:calciumpdeclosed}),
the propagating speed $\ts$ and the local curvature $\kappa$ in the curvature relation of the open-cell model~(\ref{eqn:calciumpdefull})
satisfies the system
\be\label{eqn:calciumtemp100}
\begin{split}
    (\ts + D\kappa) \frac{\di u}{\di \xi}&=D\frac{\di^2 u}{\di \xi^2}  + F(u,w) + \veps(J_{\rm in}-u),\\
    \ts \frac{\di w}{\di \xi}&=-\gamma\cdot F(u,w),
\end{split}
\ee

With the use of boundary-value solvers such as AUTO \cite{auto},
a direct numerical computation has been performed for system~(\ref{eqn:calciumtemp100})
to obtain the curvature relation of the open-cell model~(\ref{eqn:calciumpdefull})
with $J_{\rm in}=$ and $\veps=0.001$.
The result is shown in Fig.~\ref{fig:curvature_relation_open}.
Indeed,
Fig.~\ref{fig:curvature_relation_open} shows the numerically obtained curve of curvature relations superimposed on the theoretically obtained curve of curvature relations
for the closed cell model~(\ref{eqn:calciumpdeclosed}) ($\veps=0$).
It appears that in the limit $\veps\to 0$,
the curve of curvature relations of the open-cell model~(\ref{eqn:calciumpdefull}) converges to those identified in the closed-cell model.
More precisely, in the limit $\veps\to 0$,
the curve of curvature relations of the open-cell model~(\ref{eqn:calciumpdefull}) seems to arise as
the union of part of the curve $N_F$
lying above the point $\tilde{P}_m$ and the curve $N_P$.
For the stability, we can follow the argument for the closed-cell model~(\ref{eqn:calciumpdeclosed})
to obtain that the wave corresponding to the point above the point $P_T^{\rm full}$ is stable (see Fig.~\ref{fig:curvature_relation_open}).

Finally, we recall that the parabola-shaped curve $N_F$ is due to the concavity of the dispersion curve $s=s_F(\gamma)$
(see Fig.~\ref{figtemp2}), and that the existence of the curve $N_P$ is based on the existence of the wave back solutions
of the closed cell model~(\ref{eqn:calciumpdeclosed}) with wave speed $s_B(\gamma)<s_F(\gamma)$.
These features do not exist in the FHN system.

%-----------------------------------------------------------------------------------------------------

\section{Discussion and Conclusions}
\label{sec:conclusion}
\setcounter{equation}{0}

We have analyzed the dynamics of a simplified intracellular calcium (CKKOS) model. Specifically, we have investigated
the curvature relations of wave propagation in the closed-cell version~(\ref{eqn:calciumpdeclosed}) of the model  and extended the
results in conjunction with numerical computations to the curvature relations of wave propagation
for the open-cell model~(\ref{eqn:calciumpdefull}).
Comparing our results with analogous results of the canonical excitable systems, the FitzHugh-Nagumo (FHN) equations,
we find that although the temporal dynamics of the CKKOS model is  essentially the same as the FHN model,
the curvature relation of the CKKOS model does not obey the classical eikonal equation as in the FHN model.

The curvature relation is important for the evolution of waves in two spatial dimensions, e.g.,
the propagation of two-dimensional waves follows approximately the eikonal equation
in a class of excitable models including the FHN model.
The significant difference of curvature relations between the CKKONS and the FHN model suggests that
the spatio-temporal behaviours of the two models may be different.

We remark about the nature of curvature relations.
Due to the presence of the recovery variable ($\varepsilon\not=0$),
the curvature relation in the FHN system does not exactly follow the eikonal equation.
In fact, for small positive curvature parameter $\kappa$ such that $\tfrac{\epsilon}{\kappa}\ll 1$,
the results in the appendix suggests that
there are two possible propagation speed $\ts$ (see Eq.~(\ref{eqn:fhntemp115})),
although the wave corresponding to the slow speed may be unstable.
This fact has been noticed by Zykov~\cite{Zykov87}.
Further, Zykov and his coauthors~\cite{Zykov98} have shown that
for a class of the generic FHN models with large diffusivity of the recovery variable included,
the curvature relations do not follow the eikonal equation.
For such a generic FHN system with the large diffusivity of the recovery variable,
the reason for the inconsistency with the eikonal equation is not only due to the existence of the fast and slow variables,
but also to the large diffusivity of the recovery variable (see Fig.~3-4 in \cite{Zykov98}).
On the other hand, there is no recovery variable in the calcium model~(\ref{eqn:calciumpdefull}).
In fact, both $u$ and $w$ are fast variables in system~(\ref{eqn:calciumpdefull})
since the volume-ratio parameter $\gamma$ is not small.
The inconsistency with the eikonal equation for system~(\ref{eqn:calciumpdefull})
is because of the complicated dispersion relation
between wave speed $s$ and volume-ratio parameter $\gamma$ given in Fig.~\ref{figtemp2} of Sec.~\ref{sec:3.1}.
In particular, the existence of wave back solutions of system~(\ref{eqn:calciumpdeclosed})
gives rise to the existence of wave pulse solutions of system~(\ref{eqn:calciumpdeclosed}),
which in turn generates the part $N_P$ of curvature relations of system~(\ref{eqn:calciumpdeclosed}).
Such a mechanism underlying curvature relations  of system~(\ref{eqn:calciumpdeclosed})  does not exist in the generic FHN system, even with the presence of diffusivity of the recovery variable.
Finally,
earlier works on Goldbeter's model \cite{Sneyd93,Sneyd93b}, which is a model of calcium dynamics,
indicates that the curvature relation in Goldbeter's model does not follow the eikonal equation.
However, since these works are based on the piecewise linear approximation version of Goldbeter's model,
it is not clear about the mechanism constituting the inconsistency with the eikonal equation.

%%%%%%%%%%%%%%%%%%%%%%%%%%%%%%%%%%%%%%%%%%%%%%%%%%%%%%%%%%%%%%%%%%%%%%%
\bigskip
\noindent{\bf Acknowledgements.}
The authors would like to thank Professor V.S. Zykov for his valuable suggestions.
W. Zhang and J. Sneyd were supported by the Marsden Fund of the Royal Society of New Zealand.
J.-C. Tsai was partially supported by NSC and NCTS of Taiwan.

\bigskip

%%%%%%%%%%%%%%%%%%%%%%%%%%%%%%%%%%%%%%%%%%%%%%%%%%%%%%%%%%%%%%%%%%%%%%%%%%%%%

%------------------------------------------------------------------------------------------------------------
%%%%%%%%%%%%%%%%%%%%%%%%%%%%%%%%%%%%%%%%%%%%%%%%%%%%%%%%%%%%%%%%%%%%%%%%%%%%%
\appendix
\section*{Appendix: Curvature relation for the FHN system}
\def\theequation{A.\arabic{equation}}
\def\theproposition{A.\arabic{proposition}}
\setcounter{equation}{0}
\setcounter{proposition}{0}
\setcounter{section}{1}
\setcounter{subsection}{0}

\begin{figure}[!hbt]
\begin{center}
\vspace{0pt}
\includegraphics[width=250pt]{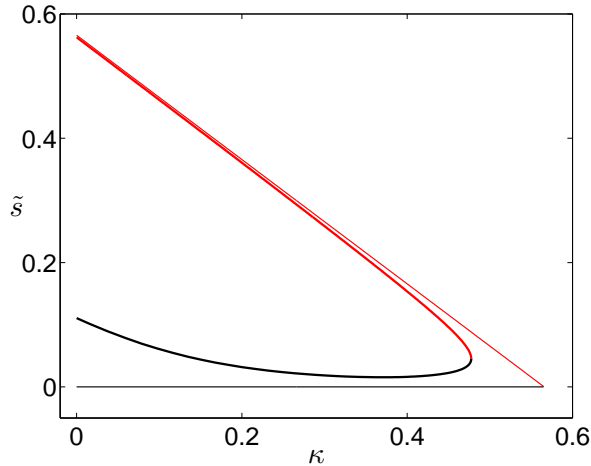}
\end{center}
\caption{\small{
Curvature relations for the FHN system~(\ref{eqn:fhnpde}),
with $D=1.0$, $\alpha=0.1$, and $\gamma=1$.
The solid curves indicate the curve of curvature relations for the case $\veps=0.01$,
while
the dashed curves shows the limit of the curve of curvature relations as $\veps\to 0$.
}}\label{fig:fhn}
\end{figure}

In this appendix, for the reader's convenience, we reproduce the curvature relation of the FHN system.
For more detail, see \cite{Zykov87}.

The classical FHN system without applied currents reads as
\be\label{eqn:fhnpde}
\begin{split}
   &\frac{\partial u}{ \partial t}=D\frac{\partial ^2 u}{\partial x^2}+f(u)-w,\\
   &\frac{\partial w}{ \partial t}=\varepsilon(u-\gamma w),
\end{split}
\ee
where
the variable $u$ represents the membrane potential, $w$ is the recovery variable.
The nonlinearity $f(u)$ is the typical bistable function, i.e.,
$f(u)=u(u-\alpha)(1-u)$ with $\alpha\in(0,\tfrac{1}{2})$.
The parameter $\varepsilon$ is a small positive number
and $\gamma$ is a positive constant.
%For numerical simulation purpose, the parameter $\alpha$ is set to be equal to 0.1 and $\gamma$ is set to be equal to 1. The %diffusion coefficient $D$ is equal to 1.

In the moving coordinate $\xi = x + s t$,
system~(\ref{eqn:fhnpde}) admits a traveling pulse $(u,w)$ with wave speed $s=s(\veps)$ such that $(u,w)(\pm\infty)=(0,0)$ and $s(\veps) = s_0  + \cal{O}(\veps)$ with $s_0 := \sqrt{D/2}\cdot(1 - 2\alpha)$ (see \cite{Jones95}).
Employing the argument as in Sec.~\ref{sec:3.2},
the relation between the propagating speed $\ts$ and the local curvature $\kappa$  is given by
\be\label{eqn:fhntemp100}
   \ts(\kappa) + D\kappa = s(\varepsilon \cdot\frac{(\ts(\kappa) + D\kappa)}{\ts(\kappa)}).
\ee
Next, we apply the Taylor expansion for the right-hand side of Eq.~\eqref{eqn:fhntemp100} around the point $\varepsilon=0$,
and then rearrange the resulting equation to get
\be\label{eqn:fhntemp105}
   (\ts(\kappa))^2 + (D\kappa  -  s_0 - s_1\varepsilon) \ts(\kappa)   - D\kappa s_1\veps = 0,
\ee
where $s_1$ is a constant depending only on the parameter $\alpha$.
Under the condition $\veps \ll \kappa \ll s_0$,
Eq.~(\ref{eqn:fhntemp105}) admits two solutions $\ts_1(\alpha)$ and $\ts_2(\alpha)$ given by
\be\label{eqn:fhntemp115}
   \ts_1(\kappa) = s_0 - D\kappa +\mathcal{O}(\veps)   \;\mbox{ and }\;  \ts_2(\kappa) = \mathcal{O}(\veps).
\ee
Up to the order $\cal{O}(\veps)$, $\ts_1(\kappa)$  is exactly the well-known eikonal equation.
Note that the wave corresponding to the speed $\ts_2(\kappa)$ is unstable (see \cite{Jones95,Zeldovich85}).
Nevertheless, the existence of the wave associated with the speed $\ts_2(\kappa)$
indicates that
the curvature relation of the FHN system~(\ref{eqn:fhnpde})
does not exactly follows the eikonal equation for small positive $\veps$.
This is consistent with numerical evidence in Fig.~\ref{fig:fhn}.

%%%%%%%%%%%%%%%%%%%%%%%%%%%%%%%%%%%%%%%%%%%%%%%%%%%%%%%%%%%%%%%%%%%%%%

\end{document}